\documentclass[12pt]{iopart}
\usepackage{iopams}
\usepackage{setstack}
\usepackage{amssymb}
\usepackage{latexsym}
\usepackage{color}
\usepackage{graphicx}   
     
\begin{document}
\title{On the Gaussian functions of two discrete variables}
\author{Nicolae Cotfas}
\address{University of Bucharest,  Physics Department,\\ P.O. Box MG-11, 077125 Bucharest, Romania}

\eads{\mailto{ncotfas@yahoo.com}}

\begin{abstract}
A remarkable discrete counterpart of the Gaussian function of one continuous  variable can be defined by using a Jacobi theta function, that is, as the sum of a convergent series.  We extend this approach to Gaussian  functions of two variables, and investigate the Fourier transform and Wigner function of the functions of discrete variable defined in this way.
\end{abstract}

\section{Introduction}

The Gaussian functions play a fundamental role in mathematics and its applications.
By using the (non-normalized) Gaussian function of continuous variable $g_\kappa :\mathbb{R}\longrightarrow \mathbb{R}$,
\begin{equation} 
g_\kappa (q)={\rm e}^{-\frac{\kappa }{2}q^2},\qquad \mbox{where $\kappa\!\in\!(0,\infty)$ is a parameter},
\end{equation}
one defines [1-5] the periodic Gaussian function (Fig. 1) of discrete variable $\mathbf{g}_\kappa :\mathbb{Z}\longrightarrow \mathbb{R}$,
\begin{equation} 
\begin{array}{l}
\mathbf{g}_{\kappa}(n)\!\!=\!\!\!\sum\limits_{\alpha =-\infty }^\infty \!\!g_\kappa \!\left( (n\!+\!\alpha d)\sqrt{ \frac{2\pi }{d}} \right)\!=\!\sum\limits_{\alpha =-\infty }^\infty{\rm e}^{-\frac{\kappa\pi}{d}(n\!+\!\alpha d)^2}
\end{array}
\end{equation}
The function $\mathbf{g}_{\kappa}$ obtained by using a method similar to Weil \cite{Weil} or Zak \cite{Zak} transform, is a generalization of Mehta's function $f_0$ \cite{Mehta}. In this article, we investigate only the case when 
$d\!=\!2j\!+\!1$ is a positive odd integer. The function $\mathbf{g}_{\kappa}$ can be written as 
\begin{equation} 
\mathbf{g}_\kappa (n)\!=\!\frac{1}{\sqrt{\kappa d}}\, \theta_3\left(\frac{n}{d},\frac{\rm i}{\kappa d} \right),
\end{equation} 
where $\theta_3$ is the Jacobi function
\begin{equation} 
\theta_3(z,\tau)\!=\!\sum\limits_{\alpha =-\infty }^\infty {\rm e}^{{\rm i}\pi \tau \alpha^2} \, {\rm e}^{2\pi {\rm i} \alpha z} 
\end{equation} 
having several remarkable properties among which we mention
\begin{equation} 
\theta_3(z,{\rm i}\tau)\!=\!\frac{1}{\sqrt{\tau}}{\rm e}^{-\frac{\pi z^2}{\tau}}\theta_3\left(\frac{z}{{\rm i}\tau},\frac{{\rm i}}{\tau}\right).
\end{equation} 
In the continuous case, the Fourier transform of $g_\kappa$ computed with the usual definition
\begin{equation} 
 \mathcal{F}[\psi](p)\stackrel{\rm def}{=}\frac{1}{\sqrt{2\pi}}\int\limits_{-\infty }^\infty {\rm e}^{-{\rm i}pq}\, \psi(q) \, dq
\end{equation} 
satisfies the relation
\begin{equation} 
\mathcal{F}[g_\kappa]\!=\!\frac{1}{\sqrt{\kappa}}\, g_{\kappa^{-1}}.
\end{equation} 
The discrete Fourier transform of $\mathbf{g}_\kappa$ , computed by using the definition 
\begin{equation} 
{\bf F}[\psi](k)\stackrel{\rm def}{=}\frac{1}{\sqrt{d}}\sum\limits_{n=-j}^j {\rm e}^{-\frac{2\pi {\rm i}}{d}kn}\, \psi(n),
\end{equation} 
satisfies a similar relation, namely \cite{Ruzzi,CD}
\begin{equation} 
\mathbf{F}[\mathbf{g}_\kappa]\!=\!\frac{1}{\sqrt{\kappa}}\, \mathbf{g}_{\kappa^{-1}}.
\end{equation} 
\begin{figure}[t]
\includegraphics[scale=0.9]{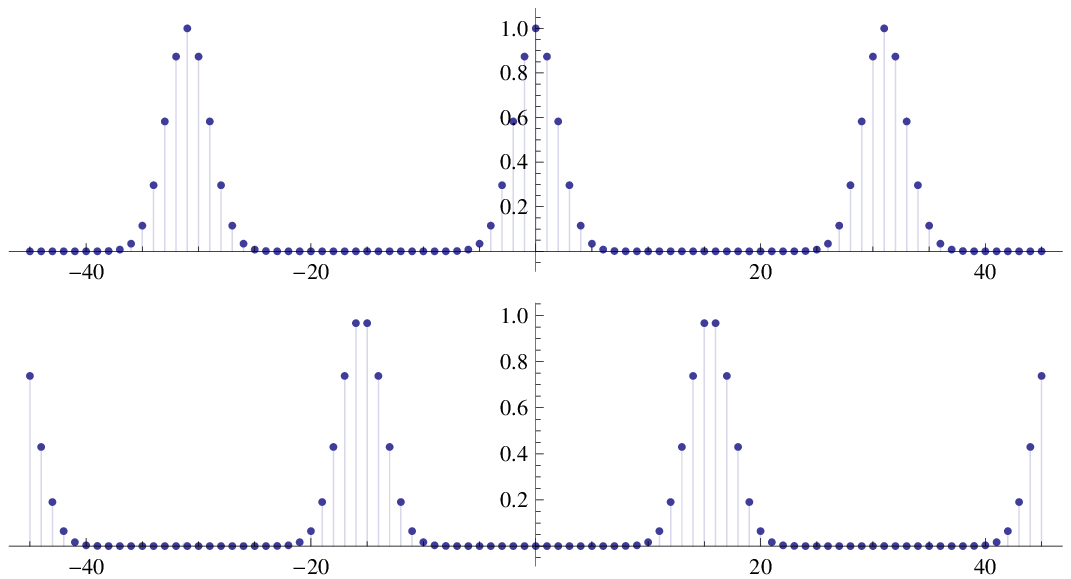} $\!\!\!\!\!\!\!\!\!\!\!\!\!\!\!\!\!\!\!\!\!\!\!\!\!\!\!\!\!\!\!\!\!\!\!\!$
\includegraphics[scale=0.5]{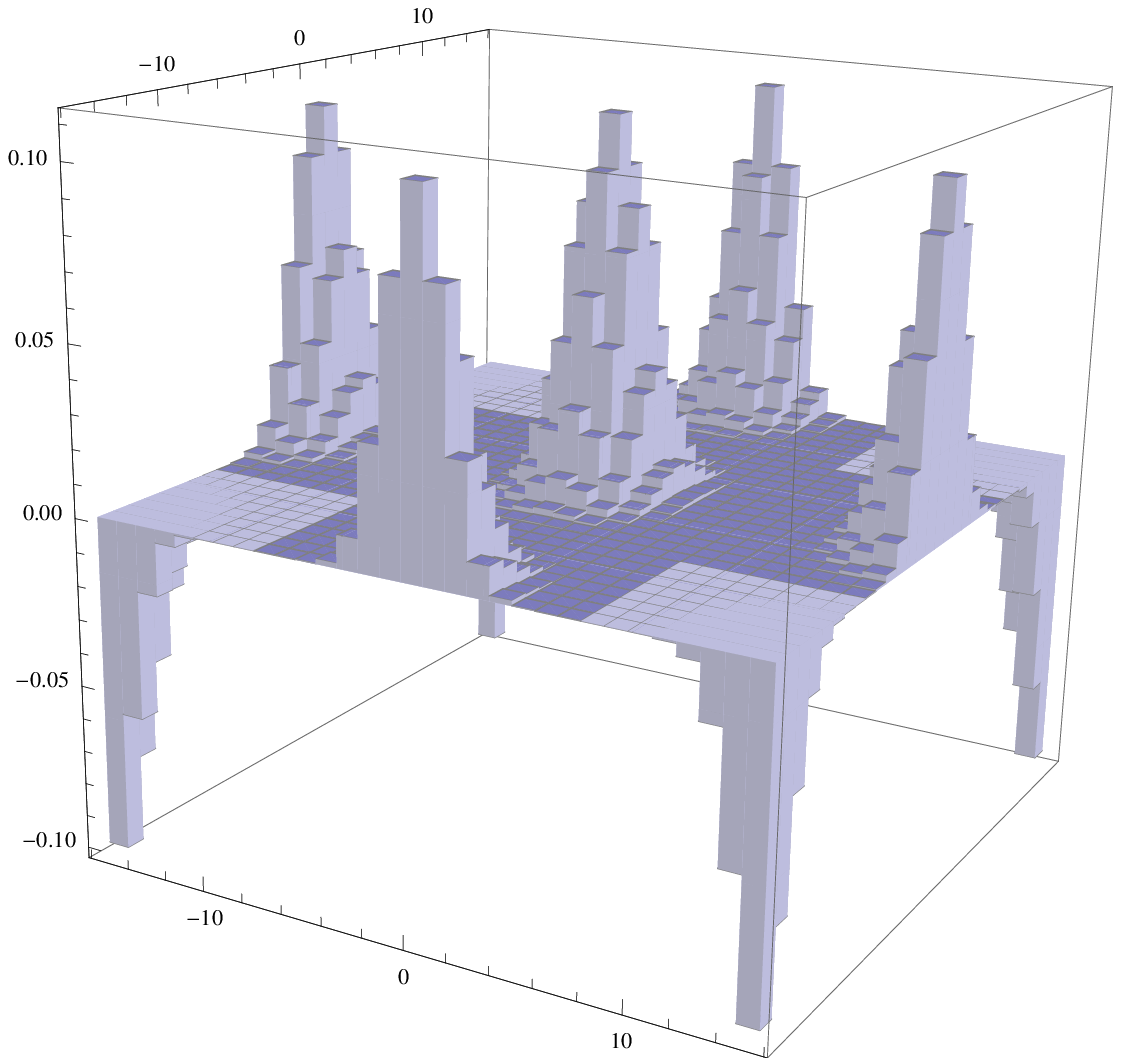} 
\caption{\label{ggk}The functions $\mathbf{g}_{\kappa}$,\,  $\mathbf{g}^+_{\kappa}$ \,  and  \, ${\bf W}_{\!{\bf g}_{\kappa} }$ \, in the case $\kappa\!=\!\frac{4}{3}$, \ $d\!=\!31$.}
\end{figure}
The Wigner function of $g_\kappa$, computed with the usual definition
\begin{equation} 
\mathcal{W}_{\psi}(q,p)\!\stackrel{\rm def}{=}\!\frac{1}{\pi }\int\limits_{-\infty}^{\infty} {\rm e}^{-2{\rm i}px}\, \overline{\psi(q\!-\!x)}\ \psi(q\!+\!x)\, dx,
\end{equation} 
is a product of Gaussian functions,
\begin{equation} 
\mathcal{W}_{g_\kappa }(q,p)\!=\!\frac{1}{\sqrt{\kappa \pi}}\, g_{2\kappa}(q)\, g_{2\kappa^{-1}}(p).
\end{equation} 
The discrete Wigner function of $\mathbf{g}_\kappa$, computed \cite{CD} by using the definition
\begin{equation} 
 {\bf W}_{\psi}(n,k)\!\stackrel{\rm def}{=}\!\frac{1}{d}\!\sum\limits_{m=-j}^j\! {\rm e}^{- \frac{4\pi {\rm i}}{d}km}\, \overline{\psi(n\!-\!m)}\ \psi(n\!+\!m)
\end{equation} 
is a sum of four products of Gaussian like functions (Fig 1)
\begin{equation} \label{Wigner1var}
\begin{array}{r}
{\bf W}_{\!{\bf g}_\kappa }\!(n,k)\!=\!\frac{1}{\sqrt{2\kappa d}}
\,\mathbf{g}_{2\kappa}(n)\, \left[\mathbf{g}_{2\kappa^{-1}}(k)\!+\!
\mathbf{g}^+_{2\kappa^{-1}}(k) \right]\ \ \\[3mm]
+\frac{1}{\sqrt{2\kappa d}}
\,\mathbf{g}^+_{2\kappa}(n)\, \left[\mathbf{g}_{2\kappa^{-1}}(k)\!-\!
\mathbf{g}^+_{2\kappa^{-1}}(k) \right],
\end{array}
\end{equation} 
where the periodic function of discrete variable $\mathbf{g}^+_\nu :\mathbb{Z}\longrightarrow \mathbb{R}$,
\begin{equation} 
\begin{array}{l}
\mathbf{g}^+_{\nu}(n)\!\!=\!\!\!\sum\limits_{\alpha =-\infty }^\infty \!\!g_\nu \!\left( (n\!+\!(\alpha\!+\!\frac{1}{2}) d)\sqrt{ \frac{2\pi }{d}} \right)\!=\!\sum\limits_{\alpha =-\infty }^\infty{\rm e}^{-\frac{\nu\pi}{d}(n\!+\!(\alpha+\frac{1}{2}) d)^2}
\end{array}
\end{equation} 
is a kind of translated Gaussian function (Fig. 1). The relation (\ref{Wigner1var}) can be written as
\begin{equation} 
\begin{array}{l}
{\bf W}_{\!{\bf g}_\kappa }\!(n,k)\!=\!C_\kappa\!\!\sum\limits_{\alpha,\beta  =-\infty }^\infty (-1)^{\alpha \beta}\ \mathcal{W}_{g_\kappa }\!\left( (n\!+\!\alpha\frac{d}{2})\sqrt{\!\frac{2\pi }{d}}, (k\!+\!\beta \frac{d}{2}) \sqrt{\!\frac{2\pi }{d}}\right),
\end{array}
\end{equation} 
where $C_\kappa$ is a constant.
Thus, there exists a simple relation between the discrete Wigner function ${\bf W}_{\!{\bf g}_\kappa }$ of a Gaussian function ${\bf g}_\kappa $ of discrete variable and the Wigner function $\mathcal{W}_{g_\kappa }$ of the corresponding Gaussian function $g_\kappa $ of continuous variable.

Our purpose is to present a version for functions of two variables of these results.

\section{Gaussian functions of two discrete variables}
Let $\sigma\!=\! {\scriptsize \left(\!\!\begin{array}{cc}
a & \!b\\
b & \!c\end{array}\!\!\right)}$ be a matrix with real entries and such that $\sigma \!>\!0$.\\
By using the Gaussian function of two continuous variables \ 
$g_\sigma\!:\!\mathbb{R}\!\times\!\mathbb{R}\!\longrightarrow \!\mathbb{R}$,
\begin{equation} 
g_\sigma(q_1,q_2)\!=\!{\rm e}^{-\frac{1}{2}(q_1\ q_2)\, {\scriptsize \left(\!\!\begin{array}{cc}
a & \!b\\
b & \!c\end{array}\!\!\right)}\,  {\scriptsize \left(\!\!\begin{array}{c}
q_1\\
q_2\end{array}\!\!\right)}}\!=\!{\rm e}^{-\frac{1}{2}[aq_1^2+2bq_1q_2+cq_2^2]}
\end{equation} 
we define the periodic {\em Gaussian function of two discrete variables} \ 
$\mathbf{g}_\sigma\!:\!\mathbb{Z}\!\times\!\mathbb{Z}\!\longrightarrow \!\mathbb{R}$,
\begin{equation} \label{gsig}
\mathbf{g}_\sigma(n_1,n_2)\!=\!\sum\limits_{\alpha_1,\alpha_2 =-\infty }^\infty 
{\rm e}^{-\frac{\pi }{d}(n_1\!+\!\alpha_1 d\ \ n_2\!+\!\alpha_2 d)\, {\scriptsize \left(\!\!\begin{array}{cc}
a & \!b\\
b & \!c\end{array}\!\!\right)}\,  {\scriptsize \left(\!\!\begin{array}{c}
n_1\!+\!\alpha_1d\\
n_2\!+\!\alpha_2 d\end{array}\!\!\right)}}
\end{equation}  
and other three complementary periodic Gaussian like functions
\begin{equation}\label{gsigp} 
\begin{array}{l}
\mathbf{g}_\sigma^{+0}(n_1,n_2)\!=\!\sum\limits_{\alpha_1,\alpha_2 =-\infty }^\infty {\rm e}^{-\frac{\pi}{d}(n_1\!+\!(\alpha_1 \!+\!\frac{1}{2})d\ \ n_2\!+\!\alpha_2 d)\, {\scriptsize \left(\!\!\begin{array}{cc}
a & \!b\\
b & \!c\end{array}\!\!\right)}\, {\scriptsize \left(\!\!\begin{array}{c}
n_1\!+\!(\alpha_1 +\frac{1}{2})d\\
n_2\!+\!\alpha_2 d\end{array}\!\!\right)}},\\[3mm] 
\mathbf{g}_\sigma^{0+}(n_1,n_2)\!=\!\sum\limits_{\alpha_1,\alpha_2 =-\infty }^\infty {\rm e}^{-\frac{\pi}{d}(n_1\!+\!\alpha_1 d\ \ n_2\!+\!(\alpha_2+\frac{1}{2}) d)\,  {\scriptsize \left(\!\!\begin{array}{cc}
a & \!b\\
b & \!c\end{array}\!\!\right)} \, {\scriptsize \left(\!\!\begin{array}{c}
n_1\!+\!\alpha_1 d\\
n_2\!+\!(\alpha_2\!+\!\frac{1}{2}) d\end{array}\!\!\right)}},\\[3mm] 
\mathbf{g}_\sigma^{++}(n_1,n_2)\!\!=\!\sum\limits_{\alpha_1,\alpha_2 =-\infty }^\infty {\rm e}^{-\frac{\pi}{d}(n_1\!+\!(\alpha_1 +\frac{1}{2})d\ \ n_2\!+\!(\alpha_2+\frac{1}{2}) d)\,  {\scriptsize \left(\!\!\begin{array}{cc}
a & \!b\\
b & \!c\end{array}\!\!\right)}\, {\scriptsize \left(\!\!\begin{array}{c}
n_1\!+\!(\alpha_1 \!+\!\frac{1}{2})d\\
n_2\!+\!(\alpha_2\!+\!\frac{1}{2}) d\end{array}\!\!\right)}}.
\end{array}
\end{equation}

\section{Discrete Fourier transform}
In the continuous case, by using the definition
\begin{equation} 
\mathcal{F}[\psi](p_1,p_2)\!\stackrel{\rm def}{=}\!\frac{1}{2\pi}\!\int\limits_{-\infty }^\infty \int\limits_{-\infty }^\infty \!{\rm e}^{-{\rm i}(p_1q_1+p_2q_2)}\, \psi(q_1,q_2) \, dq_1dq_2,
\end{equation}  
we get the known relation
\begin{equation}\label{FGauss2}
\mathcal{F}\left[{\rm e}^{-\frac{1}{2}(q_1\ \ q_2)\, \sigma \, {\scriptsize \left(\!\!\begin{array}{c}
q_1\\
q_2\end{array}\!\!\right)}} \right](p_1,p_2)\!=\!\!=\!\frac{1}{\sqrt{{\rm det\, \sigma}}}\ {\rm e}^{-\frac{1}{2}(p_1\ \ p_2)\, \sigma^{-1} \, {\scriptsize \left(\!\!\begin{array}{c}
p_1\\
p_2\end{array}\!\!\right)}},
\end{equation}
that is, we have
\begin{equation} 
\mathcal{F}[g_\sigma]\!=\!\frac{1}{\sqrt{{\rm det}\, \sigma}}\, g_{\sigma^{-1}}.
\end{equation} 
In the discrete case, the Fourier transform is usually defined as
\begin{equation} 
{\bf F}[\psi](k_1,k_2)\stackrel{\rm def}{=}\frac{1}{d}\sum\limits_{n_1=-j}^j\sum\limits_{n_2=-j}^j {\rm e}^{-\frac{2\pi {\rm i}}{d}(k_1n_1+k_2n_2)}\, \psi(n_1,n_2).
\end{equation}

\noindent {\bf Lemma 1.} {\em We have}
\begin{equation} 
\fl
\begin{array}{l}
 \mathbf{F}[\mathbf{g}_\sigma](k_1,k_2)\!=\!\frac{1}{\sqrt{{\rm det}\,\sigma}}\sum\limits_{\beta_1,\beta_2 =-\infty }^\infty \!g_{\sigma^{-1}} \!\left( (k_1\!+\!\beta_1 d)\sqrt{\frac{2\pi }{d}},(k_2\!+\!\beta_2 d)\sqrt{\frac{2\pi }{d}}  \right),\\[5mm]
\mathbf{F}[\mathbf{g}_\sigma^{+0}](k_1,k_2)\!=\!\frac{(-1)^{k_1}}{\sqrt{{\rm det}\,\sigma}}\sum\limits_{\beta_1,\beta_2 =-\infty }^\infty (-1)^{\beta_1}\ g_{\sigma^{-1}} \!\left( (k_1\!+\!\beta_1 d)\sqrt{\frac{2\pi }{d}},(k_2\!+\!\beta_2 d)\sqrt{\frac{2\pi }{d}}  \right),\\[5mm]
\mathbf{F}[\mathbf{g}_\sigma^{0+}](k_1,k_2)\!=\!\frac{(-1)^{k_2}}{\sqrt{{\rm det}\,\sigma}}\sum\limits_{\beta_1,\beta_2 =-\infty }^\infty (-1)^{\beta_2}\ g_{\sigma^{-1}} \!\left( (k_1\!+\!\beta_1 d)\sqrt{\frac{2\pi }{d}},(k_2\!+\!\beta_2 d)\sqrt{\frac{2\pi }{d}}  \right),\\[5mm]
\mathbf{F}[\mathbf{g}_\sigma^{++}](k_1,k_2)\!=\!\frac{(-1)^{k_1+k_2}}{\sqrt{{\rm det}\,\sigma}}\sum\limits_{\beta_1,\beta_2 =-\infty }^\infty (-1)^{\beta_1+\beta_2}\ g_{\sigma^{-1}} \!\left( (k_1\!+\!\beta_1 d)\sqrt{\frac{2\pi }{d}},(k_2\!+\!\beta_2 d)\sqrt{\frac{2\pi }{d}}  \right).\\[5mm]
\end{array}
\end{equation} 
{\bf Proof}.
The periodic function of two continuous variables\\[-3mm]
\begin{equation} 
\mathcal{G}_\sigma (q_1,q_2)\!=\!\sum\limits_{\alpha_1,\alpha_2 =-\infty }^\infty {\rm e}^{-\frac{\pi}{d}[a (q_1+\alpha _1d)^2+2b(q_1+\alpha_1d) (q_2+\alpha_2d)+c(q_2+\alpha_2d)^2]}
\end{equation} 
can be expended into a Fourier series\\[-3mm]
\begin{equation} 
\mathcal{G}_\sigma (q_1,q_2)\!=\!\sum\limits_{m_1,m_2=-\infty}^\infty a_{m_1m_2}\ {\rm e}^{\frac{2\pi {\rm i}}{d}(m_1q_1+m_2q_2)},\\[-3mm]
\end{equation} 
where\\[-3mm]
\[\fl 
a_{m_1m_2}\!=\!\frac{1}{d^2} \int\limits_0^d\!\! \int\limits_0^d {\rm e}^{-\frac{2\pi {\rm i}}{d}(m_1q_1+m_2q_2)}\!\!\!\sum\limits_{\alpha_1,\alpha_2 =-\infty }^\infty \!\!\!{\rm e}^{-\frac{\pi}{d}[a (q_1+\alpha _1d)^2+2b(q_1+\alpha_1d) (q_2+\alpha_2d)+c(q_2+\alpha_2d)^2]}dq_1dq_2.
\]
By using (\ref{FGauss2}) and the change of variables\ \ $q_1\!=\!y_1\sqrt{\frac{d}{2\pi}}\!-\!\alpha_1d$, \ $q_2\!=\!y_2\sqrt{\frac{d}{2\pi}}\!-\!\alpha_2d$,\ we get
\[\fl 
\begin{array}{rl}
a_{m_1m_2} & 
 \!\!\!=\!\frac{1}{2\pi d}\!\! \sum\limits_{\alpha_1,\alpha_2 =-\infty }^\infty\!\! \int\limits_{\alpha_1\sqrt{2\pi d}}^{(\alpha_1+1)\sqrt{2\pi d}}\  \int\limits_{\alpha_2\sqrt{2\pi d}}^{(\alpha_2+1)\sqrt{2\pi d}} {\rm e}^{-\frac{2\pi {\rm i}}{d}(m_1(y_1\sqrt{\frac{d}{2\pi}}\!-\!\alpha_1d)+m_2(y_2\sqrt{\frac{d}{2\pi}}\!-\!\alpha_2d))}\ {\rm e}^{-\frac{1}{2}[a y_1^2+2by_1y_2+cy_2^2]}dy_1dy_2\\[5mm]
& \!\!\!=\!\frac{1}{2\pi d}\int\limits_{-\infty}^{\infty} \int\limits_{-\infty}^{\infty} {\rm e}^{-\frac{2\pi {\rm i}}{d}(m_1y_1\sqrt{\frac{d}{2\pi}}+m_2y_2\sqrt{\frac{d}{2\pi}})}\ {\rm e}^{-\frac{1}{2}[a y_1^2+2by_1y_2+cy_2^2]}dy_1dy_2\!=\!\frac{1}{d\sqrt{{\rm det \sigma}}} {\rm e}^{-\frac{\pi}{d}(m_1\ \ m_2)\, \sigma^{-1} \, {\scriptsize \left(\!\!\begin{array}{c}
m_1\\
m_2\end{array}\!\!\right)}}.
\end{array}
\]
Consequently, we have\\[-3mm]
\begin{equation} 
\begin{array}{l}
\mathcal{G}_\sigma (q_1,q_2)\!=\!\frac{1}{d\sqrt{{\rm det\, \sigma}}}\sum\limits_{m_1,m_2=-\infty}^\infty {\rm e}^{\frac{2\pi {\rm i}}{d}(m_1q_1+m_2q_2)}\, {\rm e}^{-\frac{\pi}{d}(m_1\ \ m_2)\, \sigma^{-1} \, {\scriptsize \left(\!\!\begin{array}{c}
m_1\\
m_2\end{array}\!\!\right)}},
\end{array}
\end{equation} 
and the relation
\[\fl 
\begin{array}{rl}
\mathbf{g}_\sigma(n_1,n_2) & \!\!\!=\!\frac{1}{d\sqrt{{\rm det \, \sigma}}}\sum\limits_{m_1,m_2=-\infty}^\infty {\rm e}^{\frac{2\pi {\rm i}}{d}(m_1n_1+m_2n_2)}\, {\rm e}^{-\frac{\pi}{d}(m_1\ \ m_2)\, \sigma^{-1} \, {\scriptsize \left(\!\!\begin{array}{c}
m_1\\
m_2\end{array}\!\!\right)}}\\
& \!\!\!=\!\frac{1}{d\sqrt{{\rm det \,\sigma}}}\sum\limits_{k_1,k_2=-j}^j \sum\limits_{\beta_1,\beta_2 =-\infty }^\infty{\rm e}^{\frac{2\pi {\rm i}}{d}((k_1+\beta_1d)n_1+(k_2+\beta_2d)n_2)}\, {\rm e}^{-\frac{\pi}{d}(k_1+\beta_1d\ \ k_2+\beta_2d)\, \sigma^{-1} \, {\scriptsize \left(\!\!\begin{array}{c}
k_1+\beta _1d\\
k_2+\beta_2d\end{array}\!\!\right)}}\\
& \!\!\!=\!\frac{1}{d\sqrt{{\rm det \,\sigma}}}\sum\limits_{k_1,k_2=-j}^j {\rm e}^{\frac{2\pi {\rm i}}{d}(k_1n_1+k_2n_2)}\,\mathbf{g}_{\sigma^{-1}}(k_1,k_2)\!=\!\frac{1}{\sqrt{{\rm det \sigma}}}\mathbf{F}^{-1}[\mathbf{g}_{\sigma^{-1}}](n_1,n_2)
\end{array}
\]
equivalent to $ \mathbf{F}[\mathbf{g}_\sigma](k_1,k_2)\!=\!\frac{1}{\sqrt{{\rm det}\,\sigma}}\sum\limits_{\beta_1,\beta_2 =-\infty }^\infty \!g_{\sigma^{-1}} \!\left( (k_1\!+\!\beta_1 d)\sqrt{\frac{2\pi }{d}},(k_2\!+\!\beta_2 d)\sqrt{\frac{2\pi }{d}}  \right)$.\\
The other three relations can be obtained in a similar way. For details see \cite{NC,NC0}.\quad $\Box$\\[5mm]
\noindent {\bf Theorem 1}. {\em The discrete Fourier transform  of $\mathbf{g}_\sigma$ satisfies the relation}
\begin{equation}\label{Fgsigma}
\mathbf{F}[\mathbf{g}_\sigma]\!=\!\frac{1}{\sqrt{{\rm det}\, \sigma}}\, \mathbf{g}_{\sigma^{-1}}.
\end{equation}
{\bf Proof}. This is just the first relation from Lemma 1, written in a different way.\quad $\Box$\\[5mm]
\noindent {\bf Lemma 2.} {\em We have}
\begin{equation} \label{lemma2}
\fl 
\begin{array}{l}
\mathbf{F}[\mathbf{g}_{2\sigma}](2k_1,2k_2)\!=\!\frac{1}{2\sqrt{{\rm det}\,\sigma}}\left(
\mathbf{g}_{2\sigma^{-1}}(k_1,k_2)+\mathbf{g}_{2\sigma^{-1}}^{+0}(k_1,k_2)+
\mathbf{g}_{2\sigma^{-1}}^{0+}(k_1,k_2)+\mathbf{g}_{2\sigma^{-1}}^{++}(k_1,k_2)\right),\\[5mm]
\mathbf{F}[\mathbf{g}_{2\sigma}^{+0}](2k_1,2k_2)\!=\!\frac{1}{2\sqrt{{\rm det}\,\sigma}}\left(
\mathbf{g}_{2\sigma^{-1}}(k_1,k_2)-\mathbf{g}_{2\sigma^{-1}}^{+0}(k_1,k_2)+
\mathbf{g}_{2\sigma^{-1}}^{0+}(k_1,k_2)-\mathbf{g}_{2\sigma^{-1}}^{++}(k_1,k_2)\right),\\[5mm]
\mathbf{F}[\mathbf{g}_{2\sigma}^{0+}](2k_1,2k_2)\!=\!\frac{1}{2\sqrt{{\rm det}\,\sigma}}\left(
\mathbf{g}_{2\sigma^{-1}}(k_1,k_2)+\mathbf{g}_{2\sigma^{-1}}^{+0}(k_1,k_2)-
\mathbf{g}_{2\sigma^{-1}}^{0+}(k_1,k_2)-\mathbf{g}_{2\sigma^{-1}}^{++}(k_1,k_2)\right),\\[5mm]
\mathbf{F}[\mathbf{g}_{2\sigma}^{++}](2k_1,2k_2)\!=\!\frac{1}{2\sqrt{{\rm det}\,\sigma}}\left(
\mathbf{g}_{2\sigma^{-1}}(k_1,k_2)-\mathbf{g}_{2\sigma^{-1}}^{+0}(k_1,k_2)-
\mathbf{g}_{2\sigma^{-1}}^{0+}(k_1,k_2)+\mathbf{g}_{2\sigma^{-1}}^{++}(k_1,k_2)\right).
\end{array}
\end{equation} 
\mbox{}\\[2mm]
{\bf Proof}. The relations
\[
\fl
\begin{array}{l}
 \mathbf{F}[\mathbf{g}_{2\sigma}](2k_1,2k_2)\!=\!\frac{1}{2\sqrt{{\rm det}\,\sigma}}\sum\limits_{\beta_1,\beta_2 =-\infty }^\infty \!g_{(2\sigma)^{-1}} \!\left( (2k_1\!+\!\beta_1 d)\sqrt{\frac{2\pi }{d}},(2k_2\!+\!\beta_2 d)\sqrt{\frac{2\pi }{d}}  \right),\\[5mm]
\mathbf{F}[\mathbf{g}_{2\sigma}^{+0}](2k_1,2k_2)\!=\!\frac{1}{2\sqrt{{\rm det}\,\sigma}}\sum\limits_{\beta_1,\beta_2 =-\infty }^\infty (-1)^{\beta_1}\ g_{(2\sigma)^{-1}} \!\left( (2k_1\!+\!\beta_1 d)\sqrt{\frac{2\pi }{d}},(2k_2\!+\!\beta_2 d)\sqrt{\frac{2\pi }{d}}  \right),\\[5mm]
\mathbf{F}[\mathbf{g}_{2\sigma}^{0+}](2k_1,2k_2)\!=\!\frac{1}{2\sqrt{{\rm det}\,\sigma}}\sum\limits_{\beta_1,\beta_2 =-\infty }^\infty (-1)^{\beta_2}\ g_{(2\sigma)^{-1}} \!\left( (2k_1\!+\!\beta_1 d)\sqrt{\frac{2\pi }{d}},(2k_2\!+\!\beta_2 d)\sqrt{\frac{2\pi }{d}}  \right),\\[5mm]
\mathbf{F}[\mathbf{g}_{2\sigma}^{++}](2k_1,2k_2)\!=\!\frac{1}{2\sqrt{{\rm det}\,\sigma}}\sum\limits_{\beta_1,\beta_2 =-\infty }^\infty (-1)^{\beta_1+\beta_2}\ g_{(2\sigma)^{-1}} \!\left( (2k_1\!+\!\beta_1 d)\sqrt{\frac{2\pi }{d}},(2k_2\!+\!\beta_2 d)\sqrt{\frac{2\pi }{d}}  \right)\\[5mm]
\end{array}
\]
can be written as
\[\fl 
\begin{array}{l}
\mathbf{F}[\mathbf{g}_{2\sigma}](2k_1,2k_2)\!=\!\frac{1}{2\sqrt{{\rm det}\,\sigma}}\sum\limits_{\beta_1,\beta_2 =-\infty }^\infty \!g_{2\sigma^{-1}} \!\left( (k_1\!+\!\beta_1 \frac{d}{2})\sqrt{\frac{2\pi }{d}},(k_2\!+\!\beta_2 \frac{d}{2})\sqrt{\frac{2\pi }{d}}  \right),\\[5mm]
\mathbf{F}[\mathbf{g}_{2\sigma}^{+0}](2k_1,2k_2)\!=\!\frac{1}{2\sqrt{{\rm det}\,\sigma}}\sum\limits_{\beta_1,\beta_2 =-\infty }^\infty (-1)^{\beta_1}\ g_{2\sigma^{-1}} \!\left( (k_1\!+\!\beta_1 \frac{d}{2})\sqrt{\frac{2\pi }{d}},(k_2\!+\!\beta_2 \frac{d}{2})\sqrt{\frac{2\pi }{d}}  \right),\\[5mm]
\mathbf{F}[\mathbf{g}_{2\sigma}^{0+}](2k_1,2k_2)\!=\!\frac{1}{2\sqrt{{\rm det}\,\sigma}}\sum\limits_{\beta_1,\beta_2 =-\infty }^\infty (-1)^{\beta_2}\ g_{2\sigma^{-1}} \!\left( (k_1\!+\!\beta_1 \frac{d}{2})\sqrt{\frac{2\pi }{d}},(k_2\!+\!\beta_2 \frac{d}{2})\sqrt{\frac{2\pi }{d}}  \right),\\[5mm]
\mathbf{F}[\mathbf{g}_{2\sigma}^{++}](2k_1,2k_2)\!=\!\frac{1}{2\sqrt{{\rm det}\,\sigma}}\sum\limits_{\beta_1,\beta_2 =-\infty }^\infty (-1)^{\beta_1+\beta_2}\ g_{2\sigma^{-1}} \!\left( (k_1\!+\!\beta_1 \frac{d}{2})\sqrt{\frac{2\pi }{d}},(k_2\!+\!\beta_2 \frac{d}{2})\sqrt{\frac{2\pi }{d}}  \right).\end{array}
\]
By separating the even case from the odd  case for $\beta_1$ and  $\beta_2$, we get (\ref{lemma2}).\quad $\Box$
 \mbox{}\\[3mm]

\section{Discrete Wigner function}
The Wigner function of $g_{\sigma}$, computed by using the formula
\begin{equation} 
\fl
 \mathcal{W}_{\psi}(q_1,q_2,p_1,p_2)\!\stackrel{\rm def}{=}\!\!\frac{1}{\pi ^2}\!\int\limits_{-\infty }^\infty \int\limits_{-\infty }^\infty  {\rm e}^{-2{\rm i}(p_1x_1+p_2x_2)}\overline{\psi(q_1\!-\!x_1,q_2\!-\!x_2)} \, \psi(q_1\!+\!x_1,q_2\!+\!x_2)\, dx_1dx_2,
\end{equation} 
is a product of Gaussian functions,
\begin{equation} \label{Wignergc}
\mathcal{W}_{g_\sigma }(q_1,q_2,p_1,p_2)\!=\!\frac{1}{\pi\sqrt{{\rm det}\ \sigma}}\, g_{2\sigma}(q_1,q_2)\ g_{2\sigma^{-1}}(p_1,p_2).
\end{equation} 
 \noindent {\bf Lemma 3}. {\em If the function $f\!:\!\mathbb{R}\longrightarrow [0,\infty)$ is such that the series are convergent, then}
 \begin{equation}
 \begin{array}{l}
 \sum\limits_{\alpha,\beta=-\infty}^\infty f(\alpha,\beta)=
\sum\limits_{\mu,\eta=-\infty}^\infty f(\mu\!+\!\eta,\mu\!-\!\eta)+
\sum\limits_{\mu,\eta=-\infty}^\infty f(\mu\!+\!\eta\!+\!1,\mu\!-\!\eta).
\end{array}
\end{equation}
\mbox{}\\[2mm]
{\bf Proof}. After separating the sum as 
\begin{equation}
\sum\limits_{\alpha,\beta=-\infty}^\infty f(\alpha,\beta)=
\sum\limits_{\scriptsize \begin{array}{c}
\alpha,\beta\\
\mbox{both even}\\
\mbox{or}\\
\mbox{both odd}
\end{array}} f(\alpha,\beta)+
\sum\limits_{\scriptsize \begin{array}{c}
\alpha,\beta\\
\mbox{one even}\\
\mbox{and}\\
\mbox{other odd}
\end{array}} f(\alpha,\beta),
\end{equation}
we use the substitutions $(\alpha,\beta)\!=\!(\mu\!+\!\eta,\mu\!-\!\eta)$ and 
$(\alpha,\beta)\!=\!(\mu\!+\!\eta\!+\!1,\mu\!-\!\eta)$.\quad $\Box$\\

In the discrete case, we use for the Wigner function the definition
\begin{equation} 
\fl 
 {\bf W}_{\psi}(n_1,n_2,k_1,k_2)\!\stackrel{\rm def}{=}\!\frac{1}{d^2}\!\!\sum\limits_{m_1=-j}^j\sum\limits_{m_2=-j}^j\! \!{\rm e}^{- \frac{4\pi {\rm i}}{d}(k_1m_1+k_2m_2)}\overline{\psi(n_1\!-\!m_1,n_2\!-\!m_2)}\, \psi(n_1\!+\!m_1,n_2\!+\!m_2).
\end{equation} 
\mbox{}\\[1mm]
\noindent {\bf Theorem 2}. {\em The Wigner function ${\bf W}_{\mathbf{g}_{\sigma}}$ of $\mathbf{g}_{\sigma}$ is an algebraic sum of 16 products of Gaussian like  functions, namely}
\begin{equation}\label{Wignergd}
\fl 
\begin{array}{l}
{\bf W}_{\mathbf{g}_{\sigma}}(n_1,n_2,k_1,k_2)\\
\qquad =\frac{1}{2d\sqrt{{\rm det}\,\sigma}}\ \mathbf{g}_{2\sigma}(n_1,n_2)\left[\mathbf{g}_{2\sigma^{-1}}(k_1,k_2)+\mathbf{g}_{2\sigma^{-1}}^{+0}(k_1,k_2)+
\mathbf{g}_{2\sigma^{-1}}^{0+}(k_1,k_2)+\mathbf{g}_{2\sigma^{-1}}^{++}(k_1,k_2)\right]\\[5mm]
\qquad +\frac{1}{2d\sqrt{{\rm det}\,\sigma}}\ \mathbf{g}_{2\sigma}^{+0}(n_1,n_2)\left[
\mathbf{g}_{2\sigma^{-1}}(k_1,k_2)-\mathbf{g}_{2\sigma^{-1}}^{+0}(k_1,k_2)+
\mathbf{g}_{2\sigma^{-1}}^{0+}(k_1,k_2)-\mathbf{g}_{2\sigma^{-1}}^{++}(k_1,k_2)\right]\\[5mm]
\qquad +\frac{1}{2d\sqrt{{\rm det}\,\sigma}}\ \mathbf{g}_{2\sigma}^{0+}(n_1,n_2)\left[
\mathbf{g}_{2\sigma^{-1}}(k_1,k_2)+\mathbf{g}_{2\sigma^{-1}}^{+0}(k_1,k_2)-
\mathbf{g}_{2\sigma^{-1}}^{0+}(k_1,k_2)-\mathbf{g}_{2\sigma^{-1}}^{++}(k_1,k_2)\right]\\[5mm]
 \qquad +\frac{1}{2d\sqrt{{\rm det}\,\sigma}}\ \mathbf{g}_{2\sigma}^{++}(n_1,n_2)\left[
\mathbf{g}_{2\sigma^{-1}}(k_1,k_2)-\mathbf{g}_{2\sigma^{-1}}^{+0}(k_1,k_2)-
\mathbf{g}_{2\sigma^{-1}}^{0+}(k_1,k_2)+\mathbf{g}_{2\sigma^{-1}}^{++}(k_1,k_2)\right].\!\!
\end{array}
\end{equation}
\mbox{}\\[1mm]
{\bf Proof}. By using Lemma 3, we obtain
\[\fl
\begin{array}{l}
{\bf W}_{\mathbf{g}_{\sigma}}(n_1,n_2,k_1,k_2) =\!\frac{1}{d^2}\!\sum\limits_{m_1,m_2=-j}^j{\rm e}^{- \frac{4\pi {\rm i}}{d}(k_1m_1+k_2m_2)}\\
\qquad \qquad\qquad \qquad  \times \sum\limits_{\alpha_1,\beta_1 =-\infty }^\infty {\rm e}^{-\frac{a\pi }{d} (n_1\!-\!m_1\!+\!\alpha_1 d)^2}\, {\rm e}^{-\frac{2b\pi}{d} (n_1\!-\!m_1\!+\!\alpha_1 d)(n_2\!-\!m_2\!+\!\beta_1 d)}\, {\rm e}^{-\frac{c\pi}{d} (n_2\!-\!m_2\!+\!\beta_1 d)^2}\\
\qquad \qquad\qquad \qquad  \times \sum\limits_{\alpha_2,\beta_2 =-\infty }^\infty {\rm e}^{-\frac{a\pi }{d} (n_1\!+\!m_1\!+\!\alpha_2 d)^2}\, {\rm e}^{-\frac{2b\pi}{d} (n_1\!+\!m_1\!+\!\alpha_2 d)(n_2\!+\!m_2\!+\!\beta_2 d)}\, {\rm e}^{-\frac{c\pi}{d} (n_2\!+\!m_2\!+\!\beta_2 d)^2}\\
\qquad \quad\!=\!\frac{1}{d^2}\!\sum\limits_{m_1,m_2=-j}^j{\rm e}^{- \frac{4\pi {\rm i}}{d}(k_1m_1+k_2m_2)}\sum\limits_{\alpha_1,\alpha_2 =-\infty }^\infty\  \sum\limits_{\beta_1,\beta_2 =-\infty }^\infty{\rm e}^{-\frac{a\pi }{d} (n_1\!-\!m_1\!+\!\alpha_1 d)^2}\,  {\rm e}^{-\frac{a\pi }{d} (n_1\!+\!m_1\!+\!\alpha_2 d)^2}\\[5mm]
\qquad \qquad\times {\rm e}^{-\frac{2b\pi}{d} (n_1\!-\!m_1\!+\!\alpha_1 d)(n_2\!-\!m_2\!+\!\beta_1 d)}\, {\rm e}^{-\frac{2b\pi}{d} (n_1\!+\!m_1\!+\!\alpha_2 d)(n_2\!+\!m_2\!+\!\beta_2 d)}\ 
 {\rm e}^{-\frac{c\pi}{d} (n_2\!-\!m_2\!+\!\beta_1 d)^2} \, {\rm e}^{-\frac{c\pi}{d} (n_2\!+\!m_2\!+\!\beta_2 d)^2}
 \end{array}
\]
\[\fl
\begin{array}{r} 
\!=\!\frac{1}{d^2}\!\sum\limits_{m_1,m_2=-j}^j{\rm e}^{- \frac{4\pi {\rm i}}{d}(k_1m_1+k_2m_2)}\sum\limits_{\mu_1,\eta_1 =-\infty }^\infty\  \sum\limits_{\mu_2,\eta_2 =-\infty }^\infty{\rm e}^{-\frac{a\pi }{d} (n_1\!-\!m_1\!+\!(\mu_1+\eta_1) d)^2}\,  {\rm e}^{-\frac{a\pi }{d} (n_1\!+\!m_1\!+\!(\mu_1-\eta_1)  d)^2}\\[3mm]
{\rm e}^{-\frac{2b\pi}{d} (n_1\!-\!m_1\!+\!(\mu_1+\eta_1)  d)(n_2\!-\!m_2\!+\!(\mu_2+\eta_2)  d)}\, {\rm e}^{-\frac{2b\pi}{d} (n_1\!+\!m_1\!+\!(\mu_1-\eta_1)  d)(n_2\!+\!m_2\!+\!(\mu_2-\eta_2) d)}\\[3mm]
 {\rm e}^{-\frac{c\pi}{d} (n_2\!-\!m_2\!+\!(\mu_2+\eta_2) d)^2} \, {\rm e}^{-\frac{c\pi}{d} (n_2\!+\!m_2\!+\!(\mu_2-\eta_2) d)^2}\\
+\frac{1}{d^2}\!\sum\limits_{m_1,m_2=-j}^j{\rm e}^{- \frac{4\pi {\rm i}}{d}(k_1m_1+k_2m_2)}\sum\limits_{\mu_1,\eta_1 =-\infty }^\infty\  \sum\limits_{\mu_2,\eta_2 =-\infty }^\infty{\rm e}^{-\frac{a\pi }{d} (n_1\!-\!m_1\!+\!(\mu_1+\eta_1+1) d)^2}\,  {\rm e}^{-\frac{a\pi }{d} (n_1\!+\!m_1\!+\!(\mu_1-\eta_1)  d)^2}\\[3mm]
{\rm e}^{-\frac{2b\pi}{d} (n_1\!-\!m_1\!+\!(\mu_1+\eta_1+1)  d)(n_2\!-\!m_2\!+\!(\mu_2+\eta_2)  d)}\, {\rm e}^{-\frac{2b\pi}{d} (n_1\!+\!m_1\!+\!(\mu_1-\eta_1)  d)(n_2\!+\!m_2\!+\!(\mu_2-\eta_2) d)}\\[3mm]
 {\rm e}^{-\frac{c\pi}{d} (n_2\!-\!m_2\!+\!(\mu_2+\eta_2) d)^2} \, {\rm e}^{-\frac{c\pi}{d} (n_2\!+\!m_2\!+\!(\mu_2-\eta_2) d)^2}\\
+\frac{1}{d^2}\!\sum\limits_{m_1,m_2=-j}^j{\rm e}^{- \frac{4\pi {\rm i}}{d}(k_1m_1+k_2m_2)}\sum\limits_{\mu_1,\eta_1 =-\infty }^\infty\  \sum\limits_{\mu_2,\eta_2 =-\infty }^\infty{\rm e}^{-\frac{a\pi }{d} (n_1\!-\!m_1\!+\!(\mu_1+\eta_1) d)^2}\,  {\rm e}^{-\frac{a\pi }{d} (n_1\!+\!m_1\!+\!(\mu_1-\eta_1)  d)^2}\\[3mm]
{\rm e}^{-\frac{2b\pi}{d} (n_1\!-\!m_1\!+\!(\mu_1+\eta_1)  d)(n_2\!-\!m_2\!+\!(\mu_2+\eta_2+1)  d)}\, {\rm e}^{-\frac{2b\pi}{d} (n_1\!+\!m_1\!+\!(\mu_1-\eta_1)  d)(n_2\!+\!m_2\!+\!(\mu_2-\eta_2) d)}\\[3mm]
 {\rm e}^{-\frac{c\pi}{d} (n_2\!-\!m_2\!+\!(\mu_2+\eta_2+1) d)^2} \, {\rm e}^{-\frac{c\pi}{d} (n_2\!+\!m_2\!+\!(\mu_2-\eta_2) d)^2}\\ 
+\frac{1}{d^2}\!\sum\limits_{m_1,m_2=-j}^j{\rm e}^{- \frac{4\pi {\rm i}}{d}(k_1m_1+k_2m_2)}\sum\limits_{\mu_1,\eta_1 =-\infty }^\infty\  \sum\limits_{\mu_2,\eta_2 =-\infty }^\infty{\rm e}^{-\frac{a\pi }{d} (n_1\!-\!m_1\!+\!(\mu_1+\eta_1+1) d)^2}\,  {\rm e}^{-\frac{a\pi }{d} (n_1\!+\!m_1\!+\!(\mu_1-\eta_1)  d)^2}\\[3mm]
{\rm e}^{-\frac{2b\pi}{d} (n_1\!-\!m_1\!+\!(\mu_1+\eta_1+1)  d)(n_2\!-\!m_2\!+\!(\mu_2+\eta_2+1)  d)}\, {\rm e}^{-\frac{2b\pi}{d} (n_1\!+\!m_1\!+\!(\mu_1-\eta_1)  d)(n_2\!+\!m_2\!+\!(\mu_2-\eta_2) d)}\\[3mm]
 {\rm e}^{-\frac{c\pi}{d} (n_2\!-\!m_2\!+\!(\mu_2+\eta_2+1) d)^2} \, {\rm e}^{-\frac{c\pi}{d} (n_2\!+\!m_2\!+\!(\mu_2-\eta_2) d)^2}
 \end{array}
\]
 \[\fl
\begin{array}{r}
\!=\!\frac{1}{d^2}\!\sum\limits_{m_1,m_2=-j}^j{\rm e}^{- \frac{4\pi {\rm i}}{d}(k_1m_1+k_2m_2)}\sum\limits_{\mu_1,\eta_1 =-\infty }^\infty\  \sum\limits_{\mu_2,\eta_2 =-\infty }^\infty{\rm e}^{-\frac{2a\pi }{d} (n_1\!+\!\mu_1d)^2}\,  {\rm e}^{-\frac{2a\pi }{d} (m_1\!-\!\eta_1 d)^2}\qquad \qquad \quad \qquad  \\[3mm]
{\rm e}^{-\frac{4b\pi}{d} (n_1\!+\!\mu_1d)(n_2\!+\!\mu_2d)}\, {\rm e}^{-\frac{4b\pi}{d} (m_1\!-\!\eta_1d)(m_2\!-\!\eta_2 d)}\qquad \\[3mm]
{\rm e}^{-\frac{2c\pi }{d} (n_2\!+\!\mu_2d)^2}\,  {\rm e}^{-\frac{2c\pi }{d} (m_2\!-\!\eta_2 d)^2}\qquad \qquad \quad \qquad\\
 \!+\!\frac{1}{d^2}\!\sum\limits_{m_1,m_2=-j}^j{\rm e}^{- \frac{4\pi {\rm i}}{d}(k_1m_1+k_2m_2)}\sum\limits_{\mu_1,\eta_1 =-\infty }^\infty\  \sum\limits_{\mu_2,\eta_2 =-\infty }^\infty{\rm e}^{-\frac{2a\pi }{d} (n_1\!+\!(\mu_1\!+\!\frac{1}{2})d)^2}\,  {\rm e}^{-\frac{2a\pi }{d} (m_1\!-\!(\eta_1\!+\!\frac{1}{2}) d)^2}\quad \qquad \\[3mm]
{\rm e}^{-\frac{4b\pi}{d} (n_1\!+\!(\mu_1\!+\!\frac{1}{2})d)(n_2\!+\!\mu_2d)}\, {\rm e}^{-\frac{4b\pi}{d} (m_1\!-\!(\eta_1\!+\!\frac{1}{2}) d)(m_2\!-\!\eta_2 d)}\qquad \\[3mm]
{\rm e}^{-\frac{2c\pi }{d} (n_2\!+\!\mu_2d)^2}\,  {\rm e}^{-\frac{2c\pi }{d} (m_2\!-\!\eta_2 d)^2}\qquad \qquad \qquad\\
\!+\!\frac{1}{d^2}\!\sum\limits_{m_1,m_2=-j}^j{\rm e}^{- \frac{4\pi {\rm i}}{d}(k_1m_1+k_2m_2)}\sum\limits_{\mu_1,\eta_1 =-\infty }^\infty\  \sum\limits_{\mu_2,\eta_2 =-\infty }^\infty{\rm e}^{-\frac{2a\pi }{d} (n_1\!+\!\mu_1d)^2}\,  {\rm e}^{-\frac{2a\pi }{d} (m_1\!-\!\eta_1 d)^2}\qquad \qquad \quad\quad \\[3mm]
{\rm e}^{-\frac{4b\pi}{d} (n_1\!+\!\mu_1d)(n_2\!+\!(\mu_2\!+\!\frac{1}{2})d)}\, {\rm e}^{-\frac{4b\pi}{d} (m_1\!-\!\eta_1d)(m_2\!-\!(\eta_2\!+\!\frac{1}{2}) d)}\qquad \\[3mm]
{\rm e}^{-\frac{2c\pi }{d} (n_2\!+\!(\mu_2\!+\!\frac{1}{2})d)^2}\,  {\rm e}^{-\frac{2c\pi }{d} (m_2\!-\!(\eta_2\!+\!\frac{1}{2}) d)^2}\qquad  \quad \\
 \!+\!\frac{1}{d^2}\!\sum\limits_{m_1,m_2=-j}^j{\rm e}^{- \frac{4\pi {\rm i}}{d}(k_1m_1+k_2m_2)}\sum\limits_{\mu_1,\eta_1 =-\infty }^\infty\  \sum\limits_{\mu_2,\eta_2 =-\infty }^\infty{\rm e}^{-\frac{2a\pi }{d} (n_1\!+\!(\mu_1\!+\!\frac{1}{2})d)^2}\,  {\rm e}^{-\frac{2a\pi }{d} (m_1\!-\!(\eta_1\!+\!\frac{1}{2}) d)^2} \quad \qquad \ \ \\[3mm]
{\rm e}^{-\frac{4b\pi}{d} (n_1\!+\!(\mu_1\!+\!\frac{1}{2})d)(n_2\!+\!(\mu_2\!+\!\frac{1}{2})d)}\, {\rm e}^{-\frac{4b\pi}{d} (m_1\!-\!(\eta_1\!+\!\frac{1}{2}) d)(m_2\!-\!(\eta_2\!+\!\frac{1}{2}) d)} \\[3mm]
{\rm e}^{-\frac{2c\pi }{d} (n_2\!+\!(\mu_2\!+\!\frac{1}{2})d)^2}\,  {\rm e}^{-\frac{2c\pi }{d} (m_2\!-\!(\eta_2\!+\!\frac{1}{2}) d)^2}  \quad \qquad\ \ 
\end{array}
\]
\[\fl
\begin{array}{l}
\qquad =\frac{1}{d}{\bf g}_{2\sigma}(n_1,n_2)\ {\bf F}[{\bf g}_{2\sigma}](2k_1,2k_2)+\frac{1}{d} {\bf g}_{2\sigma}^{+0}(n_1,n_2)\ {\bf F}[{\bf g}_{2\sigma}^{+0}](2k_1,2k_2)\\[3mm]
\qquad \qquad\qquad \qquad \qquad  + \frac{1}{d}{\bf g}_{2\sigma}^{0+}(n_1,n_2)\ {\bf F}[{\bf g}_{2\sigma}^{0+}](2k_1,2k_2)+\frac{1}{d}{\bf g}_{2\sigma}^{++}(n_1,n_2)\ {\bf F}[{\bf g}_{2\sigma}^{++}](2k_1,2k_2).
\end{array}
\]
By using Lemma 2, we get (\ref{Wignergd}).\quad $\Box$

\section{Discrete-continuous correspondence}

The function $\mathbf{g}_{\sigma}$ corresponds to $g_\sigma$, namely
\begin{equation} 
\begin{array}{l}
\mathbf{g}_{\sigma}(n_1,n_2)\!=\!\sum\limits_{\alpha_1,\alpha_2 =-\infty }^\infty g_\sigma 
\left( (n_1\!+\!\alpha_1 d)\sqrt{ \frac{2\pi }{d}}, (n_2\!+\!\alpha_2 d)\sqrt{ \frac{2\pi }{d}} \right).
\end{array}
\end{equation} 
The relation (\ref{Fgsigma}) can be written as
\begin{equation} 
\begin{array}{l}
\mathbf{F}[\mathbf{g}_{\sigma}](k_1,k_2)\!=\!\sum\limits_{\beta_1,\beta_2 =-\infty }^\infty \mathcal{F}[g_\sigma] 
\left( (k_1\!+\!\beta_1 d)\sqrt{ \frac{2\pi }{d}}, (k_2\!+\!\beta_2 d)\sqrt{ \frac{2\pi }{d}} \right),
\end{array}
\end{equation} 
that is, the discrete Fourier transform $\mathbf{F}[\mathbf{g}_{\sigma}]$ of $\mathbf{g}_{\sigma}$ is the Gaussian function of two discrete variables corresponding to the continuous Fourier transform $\mathcal{F}[g_\sigma]$ of $g_\sigma$.\\
By using (\ref{gsig}), (\ref{gsigp}) and (\ref{Wignergc}), the relation (\ref{Wignergd}) can be written as
\begin{equation} 
\fl
\begin{array}{l}
{\bf W}_{\mathbf{g}_{\sigma}}(n_1,n_2,k_1,k_2)\!=\!C_\sigma \sum\limits_{\alpha_1,\alpha_2 =-\infty }^\infty \ 
\sum\limits_{\beta_1,\beta_2 =-\infty }^\infty (-1)^{\alpha_1\beta_1+\alpha_2\beta_2}\\[3mm]
\qquad \quad \qquad \qquad \quad \times\mathcal{W}_{g_{\sigma}}\left((n_1\!+\!\alpha_1\frac{d}{2})\sqrt{\!\frac{2\pi }{d}},(n_2\!+\!\alpha_2\frac{d}{2})\sqrt{\!\frac{2\pi }{d}},(k_1\!+\!\beta_1\frac{d}{2})\sqrt{\!\frac{2\pi }{d}},(k_2\!+\!\beta_2\frac{d}{2})\sqrt{\!\frac{2\pi }{d}}\right),
\end{array}
\end{equation}
where $C_\sigma$ is a constant.
Thus, the discrete Wigner function ${\bf W}_{\mathbf{g}_{\sigma}}$ of $\mathbf{g}_{\sigma}$ can be obtained directly from  the corresponding continuous Wigner function $\mathcal{W}_{g_{\sigma}}$ of $g_{\sigma}$.

\section{Concluding remarks}
Some remarkable discrete versions of the Gaussian functions, the corresponding Fourier transform  and  Wigner function can be defined  as the sum of a convergent series involving the continuous counterpart. We have investigated the Gaussian functions of two variables, but the definitions and the obtained results can easily be extended to three or more variables.

\section*{References}


\begin{thebibliography}{99}
\bibitem{Mehta}
Mehta M L 1987 Eigenvalues and eigenvectors of the finite Fourier transform {\it J. Math. Phys.} {\bf 28} 781
\bibitem{Vourdas}
Vourdas A 2004 Quantum systems with finite Hilbert space {\em Rep. Prog. Phys.} {\bf 67} 267-320
\bibitem{Ruzzi}
Ruzzi M 2006  Jacobi $\theta $-functions and discrete Fourier transform {\it J. Math. Phys.}  {\bf 47} 063507
\bibitem{CGV} 
Cotfas N, Gazeau J P and  Vourdas A 2011 Finite-dimensional Hilbert space and frame quantization {\it  J. Phys. A: Math. Theor.} {\bf  44} 175303
\bibitem{CD}
Cotfas N and Dragoman D 2012 Properties of finite Gaussians and the discrete-continuous transition {\it J. Phys. A: Math. and Theor.} {\bf 45} 425305
\bibitem{Weil} Weil A 1954 Sur certains groupes d'op\'erateurs unitaires {\it Acta math.} {\bf 111} 143-211
\bibitem{Zak}
Zak J 1967 Finite translations in solid state physics {\it Phys. Rev. Lett.} {\bf 19} 1385
\bibitem{NC}
Cotfas N 2019 Gaussian functions of two discrete variables,  \verb#https://unibuc.ro/user/nicolae.cotfas/#
\bibitem{NC0}
Cotfas N 2019 An overview on GAUSSIAN FUNCTIONS OF DISCRETE VARIABLE,  \verb#https://unibuc.ro/user/nicolae.cotfas/#



\end{thebibliography}
\end{document}